\documentclass{amsart}%
\usepackage{amsfonts}
\usepackage{amsmath}
\usepackage{amssymb}
\usepackage{graphicx}%
\newtheorem{theorem}{Theorem}
\theoremstyle{plain}

\newtheorem{lemma}{Lemma}

\theoremstyle{definition}

\newtheorem{remark}{Remark}
\theoremstyle{remark}

\numberwithin{equation}{section}

\begin{document}
\title[Asymptotic Stability of Cross Curvature Flow]
      {Asymptotic Stability of the Cross Curvature Flow at a Hyperbolic Metric}
\author{Dan Knopf and Andrea Young}
\thanks{First author partially supported by NSF grants DMS-0505920 and DMS-0545984.}

\begin{abstract}
We show that for any hyperbolic metric on a closed $3$-manifold, there exists 
a neighborhood such that every solution of a normalized cross curvature flow with 
initial data in this neighborhood exists for all time and converges to a 
constant-curvature metric. We demonstrate that the same technique proves an 
analogous result for Ricci flow. Additionally, we prove short-time existence 
and uniqueness of cross curvature flow under slightly weaker regularity 
hypotheses than was previously known.
\end{abstract}

\maketitle

It has been conjectured that any closed $3$-manifold admitting a metric with
negative sectional curvatures also admits a hyperbolic metric. This conjecture
follows from Thurston's geometrization conjecture. As is well known, Ricci
flow has shown itself to be a remarkably powerful tool in analyzing
geometrization. However, it is not expected that Ricci flow can provide a
straightforward proof of the hyperbolization conjecture. Indeed, Ricci flow
will in general preserve negative sectional curvature only in dimension $n=2$.
Moreover, one expects to see hyperbolic regions forming in a $3$-manifold only
at large times and typically after many topologically inessential surgeries.
Consequently, there are motivations to explore alternative flows that may
more readily yield information about negatively curved $3$-manifolds, and in
particular lead to a more direct proof of hyperbolization.

In 2004, Richard Hamilton and Bennett Chow proposed the cross curvature flow
(\textsc{xcf}) for $3$-manifolds  and conjectured that it would
preserve negative sectional curvature \cite{HaCh}. They further conjectured that
given any metric $g_{0}$ of strictly negative sectional curvatures on a closed
manifold $\mathcal{M}^{3}$, one would obtain a $1$-parameter family $g(t)$ of
metrics evolving by normalized \textsc{xcf} such that the metrics $g(t)$ all
have negative sectional curvatures and converge to a hyperbolic (constant
curvature) metric as $t\rightarrow\infty$. An affirmative resolution of the
Chow--Hamilton conjecture would imply that the space of hyperbolic metrics on
a fixed manifold $\mathcal{M}^{3}$ is a deformation retract of the space of
negatively-curved metrics. This conjectural picture is in sharp contrast with
the case in higher dimensions. Indeed, for $n\geq4$, Gromov and Thurston have
shown that there exist closed $n$-manifolds with sectional curvatures
$-1-\varepsilon<K\leq-1$ that admit no metric of constant curvature $K=-1$ \cite{GromThurs}.
More recently F.~Thomas Farrell and Pedro Ontaneda have shown that in all
dimensions $n\geq10$, the space of negatively curved metrics on a closed
manifold $\mathcal{M}^{n}$ has infinitely many path components whenever it is
nonempty \cite{FaOn}. Farrell and Ontaneda also provide examples of manifolds in
dimensions $n\geq10$ for which Ricci flow cannot deform all sufficiently
pinched negatively curved metrics into constant-curvature metrics \cite{FaOn1, FaOn}.

In their seminal paper \cite{HaCh}, Chow and Hamilton provide evidence in
support of their conjecture. Let $E=\operatorname*{Rc}-\frac{1}{2}Rg$ denote
the Einstein tensor and define a dual $(0,2)$-tensor $P$ by $P^{ij}%
=g^{ik}g^{j\ell}E_{k\ell}$. The functional%
\[
J(g)=\int_{\mathcal{M}^{3}}\left\{  \frac{1}{3}\operatorname*{tr}P-(\det
P)^{1/3}\right\}  \,d\mu
\]
is nonnegative by the arithmetic-geometric mean inequality and vanishes if and
only if $g$ is hyperbolic. Chow and Hamilton prove that $J(g(t))$ is
nonincreasing for as long as a solution $g(t)$ of \textsc{xcf} exists. Their
estimate is noteworthy because the effective volume $\int_{\mathcal{M}^{3}%
}\sqrt{\det P}\,d\mu$ is nondecreasing under \textsc{xcf}.

As of this writing, there are only a few other \textsc{xcf} results in the
literature. John Buckland has established short-time existence of \textsc{xcf}
for smooth initial data on compact manifolds; his approach uses De Turck
diffeomorphisms to create a strictly parabolic system \cite{Buck}.
(\textsc{xcf}, like Ricci flow, is only weakly parabolic.) Several examples of
solutions to \textsc{xcf} have been obtained by Dezhong Chen and Li Ma; these
are warped product metrics on $2$-torus and $2$-sphere bundles over the circle
\cite{MaCh}. Solutions on locally homogeneous manifolds\thinspace\footnote{In
general, \textsc{xcf} can be defined in such a way that it is weakly parabolic
if and only if the sectional curvatures of $(\mathcal{M}^{3},g)$ all have the
same sign. But on a homogeneous space, the flow reduces to an \textsc{ode}
system, so parabolicity is not an issue.} have recently been studied by
Xiaodong Cao, Yilong Ni, and Laurent Saloff-Coste \cite{CNS} and by David Glickenstein
\cite{DG}. Also, in unpublished earlier work, Ben Andrews has obtained interesting estimates
for more general solutions of \textsc{xcf}.

\bigskip

In contrast to Ricci flow, which is quasilinear, \textsc{xcf} is a fully
nonlinear system. In spite of this increased complexity, we strongly believe
that \textsc{xcf} is a highly promising tool for studying the
geometric-topological properties of negatively curved $3$-manifolds, and in
particular the hyperbolization conjecture. In this short note, we provide
further evidence in support of this belief by establishing asymptotic
stability of \textsc{xcf}. Namely, we show that for all initial data in a
sufficiently small little-H\"{o}lder $C^{2+\rho}$ neighborhood of a metric of
constant negative sectional curvature, the corresponding solution to a
suitably normalized \textsc{xcf} will exist for all time and converge to a
hyperbolic metric. (See Theorem~4 below.)\ To the best of our knowledge, these
represent the first general long-time existence and stability theorems for
\textsc{xcf}. We also prove short-time existence and uniqueness under weaker
regularity hypotheses for initial data in an appropriate neighborhood of a
hyperbolic metric. (See Theorem~3 below.) In ongoing work, we are pursuing
further applications of these results.

\bigskip

This paper is organized as follows. In \S 1, we recall the definition of
\textsc{xcf} and review some general theory regarding existence and stability
of fully nonlinear equations. In \S 2, we review little-H\"{o}lder spaces. Our
main computation is located in \S 3, where we linearize \textsc{xcf}, with a
certain normalization, about a constant curvature metric. We answer the
question of local existence and uniqueness in \S 4 and that of asymptotic
stability in \S 5. Finally, in the appendix, we apply the same methods
developed in this paper to verify asymptotic stability of Ricci flow at a
hyperbolic metric. (The stability of Ricci flow with negatively-curved initial
data satisfying certain pinching hypotheses and other geometric bounds was
studied by Rugang Ye in 1993 using alternate methods \cite{Ye}. His result is
\emph{a priori} stronger because it does not assume existence of a hyperbolic
metric. On the other hand, the approach of this paper requires no diameter or
volume hypotheses.)

\section{Fully Nonlinear Equations}
The cross curvature flow (\textsc{xcf}) is a fully nonlinear, weakly parabolic system of equations, which is defined as follows. 
As above, let $P$ denote the metric dual of the Einstein tensor. The cross curvature tensor is defined in local coordinates by
\begin{equation}
X_{ij}=\frac{1}{2}P^{k\ell}R_{ik \ell j}.
\end{equation}
Notice that if we choose an orthonormal basis so that the eigenvalues of $P$ are 
$a=-R_{2332}$, $b=-R_{1331}$, and $c=-R_{1221}$, 
then the eigenvalues of $X$ are $-bc$, $-ac$, and $-ab$. 
In this definition, our sign convention is such that $R_{ijji}$ ($i\neq j$) are the sectional curvatures, 
that is, $R_{ijk\ell}=g_{\ell m}R^m_{ijk}$. 
So if $(\mathcal{M}^{3},g)$ has negative sectional curvatures, then $X$ is negative definite. One
defines \textsc{xcf} for a manifold $(\mathcal{M}^{3},g)$ of negative sectional curvatures by
\begin{subequations}
\label{xcf}
\begin{align}
\frac{\partial g}{\partial t}&=-2X,\\
g(x,0)&=g_0(x).
\end{align}
\end{subequations}
(The opposite sign convention is used on a manifold of positive sectional curvatures.)

\bigskip

Now let us establish some notation useful for studying fully nonlinear systems in general.
Let $I$ be an interval and let $\mathbb{X}$ be a Banach space with norm $||\cdot || = ||\cdot ||_\mathbb{X}$.  For a linear operator $A:D(A) \subset \mathbb{X} \to \mathbb{X}$, we define the graph norm to be 
\[ ||x||_A = ||x|| +||Ax||, \]
\noindent 
and we let 
\[ ||A||_{L(D(A),\mathbb{X})}=\sup_{||x||_{D(A)}=1}||Ax||_\mathbb{X}. \]
\noindent
We denote the spaces of continuous and $m$ times continuously differentiable functions $f:I\to \mathbb{X}$ as $C(I;\mathbb{X})$ and $C^m(I;\mathbb{X})$ with the usual norms.  We also have weighted spaces $B_\mu((a,b];\mathbb{X})$ and $C^\alpha_\alpha((a,b];\mathbb{X})$ of functions that are bounded and H\"{o}lder continuous on $[a+\epsilon,b]$ but not necessarily up to $t=a$. 

Specifically, let $\mu \in \mathbb{R}$ and define
\[B_\mu((a,b];\mathbb{X}):= \{f:(a,b]\to \mathbb{X}: ||f||_{B_\mu((a,b];\mathbb{X})}:=\sup_{a<t\leq b}(t-a)^\mu ||f(t)|| < \infty \}.
\]

\noindent
 Let $[f]_{C^\alpha([a,b];\mathbb{X})} = \sup_{a<s<t<b}\frac{||f(t)-f(s)||}{(t-s)^\alpha}$ denote the usual $C^\alpha$ seminorm.  Then, for $0<\alpha <1$, $C^\alpha_\alpha((a,b];\mathbb{X})$ is the set of bounded functions $f:(a,b]\to \mathbb{X}$ such that 
\[ [f]_{C^\alpha_\alpha ((a,b];\mathbb{X})}:= \sup_{0<\epsilon <b-a} \epsilon^\alpha[f]_{C^\alpha([a+\epsilon,b];\mathbb{X})}< \infty,\]
with norm $||f||_{C^\alpha_\alpha ((a,b];\mathbb{X})}:= \sup_{a<t\leq b}||f(t)|| + [f]_{C^\alpha_\alpha ((a,b];\mathbb{X})}$.

We employ the theory collected in \cite{Lunardi95} regarding the local existence, uniqueness, and asymptotic behavior of solutions of fully nonlinear parabolic equations.   Let $\mathbb{D}$ be a Banach space continuously embedded in $\mathbb{X}$ and having norm $||\cdot ||_{\mathbb{D}}$.  We consider the initial value problem

\begin{equation}
\label{nonlin}
\begin{array}[c]{l}
u'(t)=F(u), \qquad{ t>0}\\
u(0)=u_0,
\end{array}
\end{equation}
\noindent
where $F: \mathcal{O} \to \mathbb{X}$ for $\mathcal{O}$ an open subset of $\mathbb{D}$.  We make several assumptions about $F$ that we will verify in \S 3 below. 

\begin{enumerate} 
\item $F$ is continuous and Fr\'{e}chet differentiable with respect to $u$.
\item The derivative $F_u$ is sectorial in $\mathbb{X}$; i.e. there are constants $\omega \in \mathbb{R}$, $\theta \in (\frac{\pi}{2}, \pi)$, $M>0$ such that $\rho(F_u) \supset S_{\theta,\omega}= \{\lambda \in \mathbb{C}: \lambda \neq \omega,|\arg(\lambda -\omega)| < \theta \}$ and 
\begin{equation}
\label{resolvent bound}
||R(\lambda, F_u)||_{L(\mathbb{X},\mathbb{X})} \leq \frac{M}{|\lambda - \omega|} 
\end{equation}
\noindent
for all $ \lambda \in S_{\theta, \omega}$.  Here $\rho (F_u)$ denotes the resolvent set of $F_u$ and $R(\lambda, F_u)=(\lambda I -F_u)^{-1}$ is the resolvent operator.
\item $F_u$ has its  graph norm  equivalent to the norm of $\mathbb{D}$.
\item Let $\bar{u} \in \mathcal{O}$.  Then there exist $r, C$ depending on $\bar{u}$ such that, for all $u,v,w \in B(\bar{u},r)$, 
\begin{equation*}
||F_u(v)-F_u(w)||_{L(\mathbb{D},\mathbb{X})}  \leq C ||v-w||_{\mathbb{D}},
\end{equation*}
where $B(\bar{u},r)$ denotes the ball around $\bar{u} \in \mathcal{O}$ of radius $r$ measured with respect to the $\mathbb{D}$ norm.   
\end{enumerate}

For such $F$, we have the following local existence and uniqueness theorem.

\begin{theorem}\cite[Theorem 8.1.1]{Lunardi95}
\label{loc exis}
Let $F(\bar{u}) \in \bar{\mathbb{D}}$.  Then there exist $\delta, r >0$, depending on $\bar{u}$, such that for all $u_0 \in B(\bar{u},r) \subset \mathbb{D}$ with $F(u_0)\in \bar{\mathbb{D}}$, there exists a solution $u$ to (\ref{nonlin}) such that $u \in C([0, \delta]; \mathbb{D})\cap C^1([0, \delta]; \mathbb{X})$.  Furthermore, $u \in C_\alpha^\alpha((0,\delta))$ and $\lim_{\epsilon \to 0}\epsilon^\alpha [u]_{C^\alpha([\epsilon, 2\epsilon];\mathbb{D})}=0$.  Finally, $u$ is the unique solution of (\ref{nonlin}) in $\bigcup_{0<\beta <1}C^\beta_\beta((0,\delta];\mathbb{D})\cap C([0,\delta];\mathbb{D})$.
\end{theorem}

We would additionally like to consider the asymptotic behavior of (\ref{nonlin}).  Notice that we can linearize this problem around a stationary solution $u$ and rewrite it as

\begin{equation}
\label{linearized}
\begin{array}[c]{l}
\bar{u}'(t)=A\bar{u}(t)+G(\bar{u}(t)-u), \qquad{ t>0}\\
\bar{u}(0)=\bar{u}_0,
\end{array}
\end{equation}
\noindent
where $A = F_{u}(\bar{u})$ and $G(\bar{u}-u)= F(\bar{u})-A\bar{u}$.  Notice that $F$ fully nonlinear means that $G$ contains ``top order'' terms.  We can assume $F(u)=0$.  We would like $A$ to be sectorial, to have graph norm equivalent to that of $\mathbb{D}$, and for the spectrum of $A$ to satisfy  
\begin{equation}
\label{spectrum}
\sup{\{\Re(\lambda): \lambda \in \sigma(A) \}} = - \omega_0 <0.
\end{equation}  

\noindent
We also want $G$ to be Fr\'{e}chet differentiable with locally Lipschitz continuous derivative and such that
\[ G( \bar{u} -u)=0, \hspace{5mm} G'( \bar{u} -u)=0. \]

\noindent
Then we have the following stability result.

\begin{theorem}
\label{stability}
Let $\omega \in [0, \omega_0)$, and let $F(\bar{u}_0) \in \bar{\mathbb{D}}$.  There exist $r, C >0$ such that for all $\bar{u}_0 \in B(u,r) \subset \mathbb{D}$ the solution $\bar{u}(t;\bar{u}_0)$ of  (\ref{nonlin}) exists for all time and 
\begin{equation}
||\bar{u}(t)-u||_{\mathbb{D}}+||\bar{u}'(t)||_{\mathbb{X}} \leq C e^{-\omega t}||\bar{u}_0||_{\mathbb{D}},
\end{equation}
\noindent
for $t \geq 0$.
\end{theorem}

\section{Little-H\"{o}lder spaces}

Let $\mathcal{M}^{3}$ denote a compact manifold admitting a hyperbolic metric
$g$. Fix a background metric $\hat{g}$ and a finite atlas $\{U_{\upsilon
}\}_{1\leq\upsilon\leq\Upsilon}$ of coordinate charts covering $\mathcal{M}%
^{3}$. For each $r\in\mathbb{N}$ and $\rho\in(0,1]$, let $\mathfrak{h}%
^{r+\rho}$ denote the little-H\"{o}lder space of symmetric $(2,0)$-tensors
with norm $\left\Vert \cdot\right\Vert _{r+\rho}$ derived from%
\[
\left\Vert u\right\Vert _{0+\rho}:=\max_{\substack{1\leq i,j\leq
3\\1\leq\upsilon\leq\Upsilon}}\left(  \sup_{x\in U_{\upsilon}}\left\vert
u_{ij}(x)\right\vert +\sup_{x,y\in U_{\upsilon}}\frac{\left\vert
u_{ij}(x)-u_{ij}(y)\right\vert }{(d_{\hat{g}}(x,y))^{\rho}}\right)  .
\]
It is well known that different choices of background metrics or atlases give equivalent norms.

Henceforth fix $\rho\in(0,1)$.  For the remainder of the paper, we will let 
\[
\mathbb{D}=\mathfrak{h}^{2+\rho}\qquad\text{and}\qquad\mathbb{X}%
=\mathfrak{h}^{0+\rho}.
\]
Then $\mathbb{D}\hookrightarrow \mathbb{X}$
is a continuous and dense inclusion.
Notice that these spaces are the closure under $||\cdot ||_{2+\rho}$ and $|| \cdot ||_{0+\rho}$ respectively 
of the space of $C^\infty$ sections of the bundle $S_2(\mathcal{M}^{3})$ of 
symmetric $(2,0)$-tensors over $\mathcal{M}^{3}$. (Recall that smooth sections are not dense
in the usual H\"older spaces.)

\section{A modified cross curvature flow}
For what follows, we consider a certain normalization of cross curvature flow which we call \textsc{kxcf}. It is defined by
\begin{equation}
\label{KXCF}
\frac{\partial \bar{g}}{\partial t}=-2X(\bar{g})-2K^2\bar{g}.
\end{equation}

\noindent
Notice that a hyperbolic metric $g$ of constant curvature $K<0$ is a fixed point of this flow.  
Such a metric is also a fixed point of the volume-normalized cross curvature flow (\textsc{nxcf}) defined to be
\begin{equation}
 \label{NXCF}
\frac{\partial \bar{g}}{\partial t}=-2X(\bar{g})+\frac{2}{3}x\bar{g},
\end{equation}
where $x=\frac{\int_{\mathcal{M}} \operatorname*{tr}_gX d\mu}{\int_{\mathcal{M}} d\mu}$.  Notice that (\ref{NXCF}) has a nonlocal term on the right-hand side.  Since both \textsc{kxcf} and \textsc{nxcf} are equivalent to
\textsc{xcf} via a reparameterization of space and time, we prefer to use the former.

\begin{lemma}
\label{rescale}
\textsc{kxcf} differs from \textsc{xcf} only by a change of scale in space and time.  
\end{lemma}
\begin{proof}
Define dilating factors $\psi (t)>0$ by $\psi(t)=Ae^{2K^2t}$ and 
define $\tilde{t}=\int_0^t \psi^2(\tau)\,d\tau$, so that $\frac{d\tilde{t}}{dt}=\psi^2(t)$.  
If we let $\tilde{g}=\psi \bar{g}$, then $X(\tilde{g})=\frac{1}{\psi}X(\bar{g})$.  
Supposing $\bar{g}$ solves (\ref{KXCF}), we have the following computation:

\begin{eqnarray*}
\frac{\partial \tilde{g}}{\partial \tilde{t}} = \frac{dt}{d \tilde{t}}(\frac{\partial}{\partial t}(\psi \bar{g}))&=& \frac{dt}{d \tilde{t}}(\frac{\partial \psi}{\partial t}\bar{g})+ \frac{dt}{d \tilde{t}}(\psi(-2X(\bar{g})-2K^2\bar{g}))\\
&=&\frac{1}{\psi^3}\frac{\partial \psi}{\partial t}\tilde{g}-2X(\tilde{g})-\frac{2}{\psi^2}K^2\tilde{g}\\
&=&-2X(\tilde{g}).
\end{eqnarray*}
\noindent
Thus $\tilde{g}$ solves (\ref{xcf}), and we have shown the desired equivalence.
\end{proof}

We now want to define a DeTurck-modified cross curvature flow%
\[
\frac{\partial}{\partial t}\bar{g}(x,t)=F(x,\bar{g}(x,t))
\]
for Riemannian metrics $\bar{g}(\cdot,t)$ in a neighborhood $\mathcal{O}%
\subset\mathbb{D}$ of the hyperbolic metric $g$ on $\mathcal{M}^{3}$. Here
$\mathcal{O}$ is an open set in $\mathbb{D}$ to be determined below.

Given $\bar{g}\in\mathcal{O}$ and a smooth section $h$ of $S_{2}(\mathcal{M}^{3})$,
define a vector field $Y(\bar{g},h)$ on $\mathcal{M}^{3}$ in local coordinates
by%
\begin{equation}
Y^{\ell}(\bar{g},h):=\frac{1}{2}\bar{g}^{k\ell}\partial_{k}(\bar{g}^{ij}%
h_{ij})-\bar{g}^{k\ell}\bar{g}^{ij}\bar{\nabla}_{i}h_{jk}.
\end{equation}
Assume that $g$ has constant sectional curvature $K<0$, and consider the
\emph{DeTurck cross curvature flow} (\textsc{dxcf}) given by
\begin{subequations}
\label{XXCF}%
\begin{align}
\frac{\partial}{\partial t}\bar{g} &  =F(\bar{g}):=-2X(\bar{g})+K\mathcal{L}%
_{Y(g,\bar{g})}g-2K^{2}\bar{g}\\
\bar{g}(0) &  =g_{0}.
\end{align}
Notice that $F(g)=0$.

We next derive the linearization of \textsc{dxcf}.
\end{subequations}

\begin{lemma}
If $\tilde{g}=\bar{g}+h$, the Fr\'{e}chet derivative $F_{\bar{g}}$ is the
linear operator $A_{\bar{g}}$ given by
\begin{align*}
(A_{\bar{g}}h)_{ik} &  =\frac{1}{2}\bar{R}_{jik}^{\ell}\{\bar{\Delta}_{\ell
}h_{\ell}^{j}+(\mathcal{L}_{Y(\bar{g},h)}\bar{g})_{\ell}^{j}\}-\frac{1}{4}%
\bar{R}\{\bar{\Delta}_{\ell}h_{ik}+(\mathcal{L}_{Y(\bar{g},h)}\bar{g}%
)_{ik}\}\\
&  +\frac{1}{2}\bar{R}_{\ell}^{j}(\bar{\nabla}_{i}\bar{\nabla}_{k}h_{j}^{\ell
}-\bar{\nabla}_{j}\bar{\nabla}_{k}h_{i}^{\ell}-\bar{\nabla}_{i}\bar{\nabla
}^{\ell}h_{jk}+\bar{\nabla}_{j}\bar{\nabla}^{\ell}h_{ik})\\
&  -\frac{1}{2}\bar{R}_{ik}(\bar{\Delta}H-\bar{\delta}^{2}h)+K(\mathcal{L}%
_{Y(\bar{g},h)}\bar{g})_{ik}-2K^{2}h\\
&  -\bar{R}_{ijk}^{\ell}\bar{R}_{\ell}^{m}h_{m}^{j}+\frac{1}{2}\bar{R}%
_{ijm}^{\ell}\bar{R}_{\ell}^{j}h_{k}^{m}-\frac{1}{2}\bar{R}_{ijk}^{m}\bar
{R}_{\ell}^{j}h_{m}^{\ell}-\frac{1}{2}\left\langle \bar
{\operatorname*{Rc}},h\right\rangle _{\bar{g}}\bar{R}_{ik}.
\end{align*}
\noindent
At a metric $g$ of constant sectional curvature $K<0$, one has
\begin{center} 
$R_{ijk\ell}=K(g_{i\ell}g_{jk}-g_{ik}g_{j\ell})$, $\operatorname*{Rc}=2Kg$, and $R=6K$,
\end{center} 
whence the formula above reduces to
\[
A_{g}h=-K\Delta h-2K^{2}Hg+2K^{2}h,
\]
where $H=g^{ij}h_{ij}$.
\end{lemma}
\begin{proof}
Using the definition of the cross curvature tensor, $X_{ij}=\frac{1}{2}P^{k\ell}R_{ik\ell j}$, one linearizes using standard first variation formulas, as found for instance in \cite{Besse} or \cite{GIK}.  
\end{proof}

Observe that $A_g$ is a self-adjoint elliptic operator. The $L^{2}$ spectrum of $A_{g}$ consists
of discrete eigenvalues of finite multiplicity contained in the half-line
$(-\infty,2K^{2}]\ $and accumulating only at $-\infty$. Standard Schauder
theory implies that $A_{g}$ is sectorial with its graph norm equivalent to $\left\Vert \cdot\right\Vert _{2+\rho}$.  In particular,
there exists $C\in(0,\infty)$ such that%
\begin{equation}
\frac{1}{C}\left\Vert h\right\Vert _{A_{g}}\leq\left\Vert h\right\Vert
_{2+\rho}\leq C\left\Vert h\right\Vert _{A_{g}}.\label{graphnorm}%
\end{equation}

Noting that $\Delta_{\ell}h=\Delta h+H\operatorname*{Rc}-Rh$ on $(\mathcal{M}%
^{3},g)$, one may also write $A_{g}$ in the form%
\[
A_{g}h=-K(\Delta_{\ell}h+4Kh).
\]

\section{Local existence and uniqueness}

Recall our notation that $\mathbb{D}=\mathfrak{h}^{2+\rho}$ and $\mathbb{X}%
=\mathfrak{h}^{0+\rho}$.  Let $\bar{g}\in\mathbb{D}$. In each coordinate chart $U_{\upsilon}$, one may
write%
\begin{equation}
(A_{\bar{g}}h)_{ij}=a^{k\ell}\partial_{k}\partial_{\ell}h_{ij}+b^{k}%
\partial_{k}h_{ij}+c_{ij}^{k\ell}h_{k\ell}, \label{local}%
\end{equation}
where $a$, $b$, and $c$ depend on $x\in U_{\upsilon}$ and $\bar
{g},\partial\bar{g},\partial^{2}\bar{g}$. By taking $\bar{g}$ close enough to
$g$ in $\mathbb{D}$, we can make $a,b,c$ as close in $L^{\infty}$ as desired
to their values for $A_{g}$.

Let $\mathcal{O}:=B_{\eta}^{2+\rho}(g)$, where $B_{\eta}^{2+\rho}(g)$ denotes the
$\eta$-ball around $g$; i.e.
\[
B_{\eta}^{2+\rho}(g):=\{\bar{g}\in\mathfrak{h}^{2+\rho}:\left\Vert \bar
{g}-g\right\Vert _{2+\rho}<\eta\}. 
\] 
Choose $\eta>0$ small enough such
that for all $\bar{g}\in\mathcal{O}$,

\begin{enumerate}
\item $\bar{g}$ is a Riemannian metric,

\item $A_{\bar{g}}$ is uniformly elliptic, and

\item there exists a sufficiently small $\delta >0$, to be chosen below, such that $\left\Vert (A_{\bar{g}}-A_{g})h\right\Vert _{0+\rho}<\delta
\left\Vert h\right\Vert _{2+\rho}$.  
\end{enumerate}

Let $\delta=(M+1)^{-1}$, with $M$ as in (\ref{resolvent bound}) depending only on the maximum of the resolvent operator.  Then it is a standard fact that $A_{\bar{g}}$ is sectorial for all $\bar{g}%
\in\mathcal{O}$. (For example, see \cite[Proposition~2.4.2]{Lunardi95}.)  We can then choose $\delta$ smaller if necessary (depending on $C$ in ({\ref{graphnorm})) so that the graph norm of $A_{\bar{g}}$ is equivalent to $||\cdot||_{2+\rho}$.

Let
\[ G(h)=F(g+h)-A_gh= -2X(g+h)+K\mathcal{L}_{Y(g,g+h)}g-2K^2(g+h)-A_gh. \]  Since g is a fixed point of (\ref{KXCF}), and $A_g$ is linear, we see that  $G(0)=0$.  From the computation above, it is clear that the derivative of $G$ at $h$ is  
\[ G_hk=A_{g+h}k-A_gk,\]
so $G'(0)=0$ as well.  The fact that for any $r \in (0,\eta]$, there exists $C>0$ such that $||G_hz||_\mathbb{X}\leq C||z||_\mathbb{D}$ uniformly for $h \in B_r^{2+\rho}(0)$ follows from property~(3) above.
This establishes the local Lipschitz continuity that we need to apply Theorem~\ref{stability}.

Given $\bar{g}\in\mathcal{O}$, choose $\varepsilon>0$ small enough that
$B_{\varepsilon}^{2+\rho}(\bar{g})\subseteq\mathcal{O}$. Fix any coordinate
chart $U_{\upsilon}$. Given $u\in B_{\varepsilon}^{2+\rho}(\bar{g})$, let
$a(x)\equiv a(x,u,\partial u,\partial^{2}u)$, $b(x)\equiv b(x,u,\partial
u,\partial^{2}u)$, and $c(x)\equiv c(x,u,\partial u,\partial^{2}u)$ denote the
local coefficients of $A_{u}$, as in (\ref{local}). Then for any
$h\in\mathbb{D}$, one has%
\begin{align*}
\frac{\left\vert a^{k\ell}(x)\partial_{k}\partial_{\ell}h_{ij}(x)-a^{k\ell
}(y)\partial_{k}\partial_{\ell}h_{ij}(y)\right\vert }{(d_{\hat{g}}%
(x,y))^{\rho}} &  \leq\left\vert a^{k\ell}(x)\frac{\partial_{k}\partial_{\ell
}h_{ij}(x)-\partial_{k}\partial_{\ell}h_{ij}(y)}{(d_{\hat{g}}(x,y))^{\rho}%
}\right\vert \\
&  +\left\vert \partial_{k}\partial_{\ell}h_{ij}(y)\frac{a^{k\ell}%
(x)-a^{k\ell}(y)}{(d_{\hat{g}}(x,y))^{\rho}}\right\vert \\
&  \leq\left\Vert u\right\Vert _{2+0}\left\Vert h\right\Vert _{2+\rho
}+\left\Vert h\right\Vert _{2+0}\left\Vert u\right\Vert _{2+\rho}.
\end{align*}
Similar arguments apply to the lower order terms.  For the details, one need only make minor modifications to the proof of Lemma 3.3 in \cite{GIK}.  In this way, it is easy to see that%
\[
\left\Vert A_{u}(v)-A_{u}(w)\right\Vert _{0+\rho}\leq C\left\Vert u\right\Vert
_{2+\rho}\left\Vert v-w\right\Vert _{2+\rho}%
\]
for all $u,v,w\in B_{\varepsilon}^{2+\rho}(\bar{g})$. 

Then we can apply
Theorem \ref{loc exis} to obtain the following theorem.

\begin{theorem}
Let $(\mathcal{M}^{3},g)$ be a Riemannian manifold having constant sectional curvature $K <0$.  
There exist $\delta, r >0$ such that for all $\bar{g}_0 \in B_r^{2+\rho}(g)$ there exists a 
solution $\bar{g} \in C([0, \delta]; \mathfrak{h}^{2+\rho})\cap C^1([0,\delta];\mathfrak{h}^{0+\rho})$ 
of \textsc{dxcf} that exists for all $t\in [0,\delta]$.  This is the unique solution 
in $\bigcup_{0<\beta <1}C^\beta_\beta((t_0,t_0+\delta];\mathfrak{h^{2+\rho}})\cap C([t_0,t_0+\delta];\mathfrak{h^{2+\rho}})$.
\end{theorem}

\noindent
\begin{remark}
The idea of using a DeTurck trick to prove short-time existence and uniqueness
for \textsc{xcf} is due to Buckland \cite{Buck}. Because of the invariance of
the Riemann curvature tensor under the infinite-dimensional diffeomorphism
group, the symbol of the linearization of (\ref{xcf}) is only degenerate
elliptic. Indeed, Buckland shows that, in appropriate coordinates, the symbol
of the linearization, acting on a vector $(h_{11},h_{12},h_{13},h_{22}%
,h_{33},h_{23})^{T}$ representing the components of a variation, is
represented by the matrix%
\[%
\begin{pmatrix}
0 & 0 & 0 & \Lambda^{22} & \Lambda^{33} & 2\Lambda^{23}\\
0 & 0 & 0 & -\Lambda^{12} & 0 & -\Lambda^{13}\\
0 & 0 & 0 & 0 & -\Lambda^{13} & -\Lambda^{12}\\
0 & 0 & 0 & \Lambda^{11} & 0 & 0\\
0 & 0 & 0 & 0 & \Lambda^{11} & 0\\
0 & 0 & 0 & 0 & 0 & \Lambda^{11}%
\end{pmatrix}
\]
with eigenvalues $\Lambda^{11}>0$ and $0$, each of multiplicity three. See
\cite[equation~(4)]{Buck}. Ellipticity fails because of the null eigenvalue.
The DeTurck diffeomorphisms effectively fix a gauge, breaking diffeomorphism
invariance and making the linearization strongly elliptic. In this way,
Buckland proves short-time existence and uniqueness  for smooth initial data.
The theorem above may be regarded as a mild extension of his result, in the
sense that it proves existence and uniqueness of solutions to \textsc{dxcf}
for somewhat less regular initial data, at least for such data sufficiently
near a hyperbolic metric.
\end{remark}

\section{Stability}

Without loss of generality, we may assume that $(\mathcal{M}^{3},g)$ has
constant sectional curvature $K=-1$. Henceforth write $A\equiv A_{g}$, noting
that%
\begin{align*}
Ah &  =\Delta_{\ell}h-4h\\
&  =\Delta h-2Hg+2h.
\end{align*}
Clearly, the $L^{2}$ spectrum of $A$ is contained in $(-\infty,\omega_{0}]$
for some $\omega_{0}\leq2$.  We now further analyze the spectrum, using a Bochner formula due to Koiso \cite{Koiso}.  
Notice that, for $h$ a symmetric $(2,0)$-tensor on a closed manifold $(\mathcal{M}^{n},g)$, one has
\[||\nabla h||^2=||\delta h||^2+\frac{1}{2}||T||^2+\int_{\mathcal{M}^{n}}(R_{ijkl}h^{il}h^{jk}-R_i^kh_{jk}h^{ij})d\mu, \]
where $T=T(h)$ is a $(3,0)$-tensor defined by $T_{ijk}=\nabla_k h_{ij}-\nabla_i h_{jk}$ and $(\delta h)_k =-g^{ij}\nabla_i h_{jk}$.  In our case, $g$ has constant sectional curvature $K<0$, so again we use the formulas $R_{ijk\ell}=K(g_{i\ell}g_{jk}-g_{ik}g_{j\ell})$ and $\operatorname*{Rc}=2Kg$.  Thus Koiso's Bochner formula reduces to
\[|| \nabla h||^2 = || \delta h||^2+\frac{1}{2}||T||^2-||H||^2+3||h||^2. \]
\noindent
This observation implies that
\[ \int (Ah,h) d\mu \leq -||H||^2-||h||^2 \leq -||h||^2<0. \]
\noindent
Thus there exists an $\omega\geq 1$ such that the $L^2$ spectrum of $A_g$ is contained in the half-line $(-\infty, -\omega]$.  So we can apply Theorem \ref{stability} to obtain 
asymptotic stability for \textsc{dxcf}.

Finally, we shall show that asymptotic stability for \textsc{dxcf} implies the same for \textsc{xcf}.  
We will utilize the following lemma, whose proof may be found in \cite{GIK}.

\begin{lemma}
\label{dansoldlemma}
Let $Y(t)$ be a vector field on a Riemannian manifold $(\mathcal{M}^{n},g(t))$, where $0\leq t<\infty$, and suppose there are constants $0<c\leq C <\infty$ such that
\[\sup_{x\in M^n}|Y(x,t)|_{g(t)}\leq Ce^{-ct}.
\]
\noindent
Then the diffeomorphisms $\varphi_t$ generated by $Y$ converge exponentially to a fixed diffeomorphism $\varphi_\infty$ of $\mathcal{M}^{n}$.
\end{lemma}  

\begin{lemma}
Let $g$ be metric of constant negative sectional curvature on $\mathcal{M}^{3}$.  Suppose there exists an $r$ such that for all $\tilde{g}_0 \in B_r^{2+\rho}(g)$, the unique solution $\bar{g}(t)$ of (\ref{XXCF}) with $\bar{g}(0)=\tilde{g}_0$ converges exponentially fast to $g$.  Then the unique solution $\tilde{g}(t):=\varphi_t^*\bar{g}$ of  (\ref{KXCF}) with $\tilde{g}(0)=\tilde{g}_0$ converges exponentially fast to a constant curvature metric $\tilde{g}_\infty$.
\end{lemma}
\begin{proof}
Recall that $Y$ is defined to be
\[Y^l=\frac{1}{2}\bar{g}^{kl}\partial_k(g^{ij}\bar{g}_{ij})-g^{kl}g^{ij}\nabla_i\bar{g}_{jk}.\]
\noindent
Since $\bar{g}(t)\to g$ exponentially fast, we have $Y^l\to 0$ exponentially fast as well.
So Lemma \ref{dansoldlemma} implies that the diffeomorphisms $\varphi_t$ converge to a fixed diffeomorphism $\varphi_\infty$.  Thus $\tilde{g}(t)$ converges to a limit metric $\tilde{g}_\infty$, which by diffeomorphism invariance has constant sectional curvature.
\end{proof}

Then we can apply Theorem \ref{stability} to obtain asymptotic stability of \textsc{kxcf}.

\begin{theorem}
Let $(\mathcal{M}^{3},g)$ be a closed Riemannian manifold with constant sectional curvature $K<0$.  
Then there exists $\delta>0$ such that for all $\bar{g}_0\in B_\delta^{2+\rho}(g)$, the solution $\bar{g}$ to (\ref{KXCF}) having initial condition $\bar{g}_0$ exists for all time and converges exponentially fast to a constant curvature hyperbolic metric.
\end{theorem}

\appendix
\section{Asymptotic Stability of Ricci flow at a hyperbolic metric}

The methods developed in this paper provide a simple proof of the asymptotic stability of Ricci flow at a hyperbolic metric.  As noted above, a more powerful stability result was obtained earlier by Ye using somewhat different methods \cite{Ye}. Recall that the Ricci flow is defined to be
\begin{subequations}
\label{ricci flow}
\begin{align}
\frac{\partial \bar{g}}{\partial t}&=-2\operatorname{Rc}(\bar{g}), \\
\bar{g}(0)&=\bar{g}_0.
\end{align}
\end{subequations}
\noindent
We proceed as above and define a normalized Ricci flow (\textsc{knrf}) that differs from the usual volume-normalized 
flow but which also can be obtained from Ricci flow only by a reparameterization of space and time.  \textsc{knrf} is given by

\begin{subequations}
\label{krf}
\begin{align}
\frac{\partial \bar{g}}{\partial t} &=-2\operatorname{Rc}(\bar{g})+4K\bar{g},\\
\bar{g}(0) &= \bar{g}_0.
\end{align}
\end{subequations}
\noindent
In particular, a constant curvature metric with $K<0$ is a fixed point of \textsc{knrf}.

Using standard variation formulas, we can linearize the right hand side of the \textsc{knrf} equation 
about a hyperbolic metric $g$ having constant curvature $K=-1$ to obtain
\begin{equation}
\label{ricci lin}
A_gh= \Delta_\ell h +4Kh = \Delta h -2Hg +2h.
\end{equation}

The same trick as above allows us to bound the spectrum in the interval $(\infty, -1]$.  By mimicking our previous analysis, one easily checks that the hypotheses of Theorem \ref{stability} are satisfied.
(Or see the detailed calculations in \cite{GIK}, which treat a more technically difficult case in which there is a center manifold present.) Thus we obtain the following theorem.

\begin{theorem}
Let $(\mathcal{M}^{3},g)$ be a closed Riemannian manifold having constant sectional curvature $K<0$.
Then there exists $\delta>0$ such that for all $\bar{g}_0 \in B_\delta^{2+\rho}(g)$, the solution $\bar{g}$ to (\ref{krf}) having initial condition $\bar{g}_0$ exists for all time and converges exponentially fast to a constant curvature hyperbolic metric.
\end{theorem}
Notice that curvature pinching is not quite sufficient to apply this theorem: one also requires a H\"{o}lder bound on the 
first derivative of the curvature. Compare \cite{FaOn1}.

\end{document}